 \documentclass[12pt,a4paper,twoside]{article}

\usepackage{amsfonts,amsmath,amssymb,amsthm,latexsym}
\usepackage{times,mathptm}
\usepackage{latexsym}
\usepackage[all]{xy}
\usepackage{amsfonts}
\usepackage{amsthm}
\usepackage{amsmath}
\usepackage{amssymb}
\usepackage{amscd}
 \numberwithin{equation}{section}

  \def\sw#1{{\sb{(#1)}}} 

  \def\<{{\langle}} 
  \def\>{{\rangle}} 
   
  \def\eps{\varepsilon}

  \def\note#1{{}}

  \def\note#1{} 
  \def\M{{\bf M}}

  \def\cC{{\mathcal C}} 
  \def\cA{{\mathcal A}} 
  \def\bA{\bar{\cA}}
  \def\ba{\bar{a}}

  \def\beq{\begin{equation}} 
  \def\eeq{\end{equation}} 
  \def\DC{{\Delta_\cC}} 
  \def \eC{{\eps_\cC}}

  \def\im{{\rm Im}}

  \def\ut{{\otimes}} 
  \def\ot{{\otimes}}

  \def\roM{\varrho^{M}} 
\def\Nro{{}^{N}\!\varrho}
\def\Mro{{}^{M}\!\varrho} 
  \def\roL{\varrho^{L}} 
\def\Lro{{}^{L}\varrho}

  \def\roA{\varrho^{A}}
 \def\Aro{{}^{A}\varrho}
  \def\robA{\varrho^{\bA}}
 \def\bAro{{}^{\bA}\varrho}
 \def\Vect{\mathbf{Vect}}
 
 \def\coker{\mathrm{coker}}
 \def\eA{{\eta_\cA}}
    \def\mA{{\mu_\cA}}
    \def\roA{{\varrho^\cA}}

    \def\Aro{{}^\cA\!\varrho}
    \def\coten#1{\square_{#1}}
    \def\wra{\varrho_M}

  \newcounter{zlist} 
  \newenvironment{zlist}{\begin{list}{(\arabic{zlist})}{ 
  \usecounter{zlist}\leftmargin2.5em\labelwidth2em\labelsep0.5em 
  \topsep0.6ex
  \parsep0.3ex plus0.2ex minus0.1ex}}{\end{list}}

  \newcounter{blist}

  \newcounter{rlist}

\numberwithin{equation}{section}

  \def\Label#1{\label{#1}\ifmmode\llap{[#1] }\else 
  \marginpar{\smash{\hbox{\tiny [#1]}}}\fi} 
  \def\Label{\label}

  \newcounter{c} 
   
  \newcommand{\etyk}[1]{\vspace{-7.4mm}$$\begin{equation}\Label{#1} 
  \addtocounter{c}{1}} 
  \renewcommand{\]}{\ifnum \value{c}=1 $$\else \end{equation}\fi} 
\newtheorem{proposition}{Proposition}
\newtheorem{lemma}[proposition]{Lemma}

\newtheorem{theorem}[proposition]{Theorem}

\theoremstyle{definition}
\newtheorem{definition}[proposition]{Definition}

\theoremstyle{remark}

\renewcommand{\thefootnote}{\fnsymbol{footnote}}
\title{\MakeUppercase{Flat connections and (co)modules}}
\author{Tomasz Brzezi\'nski}
\date{}
\renewcommand{\date}{\vspace{-5mm}}
\begin{document}
\maketitle \vspace*{-3mm}\relax
\renewcommand{\thefootnote}{\arabic{footnote}}

\noindent \textit{\small Department of Mathematics, University of Wales Swansea,
   Singleton Park, \newline\indent  Swansea SA2 8PP, U.K.\\
e-mail: T.Brzezinski@swansea.ac.uk}

\begin{abstract}
 The relationship between comodules of a coring and flat connections is reviewed. In particular we specialise to corings which are built on a tensor product of algebra and a coalgebra. Such corings are in one-to-one correspondence with entwining structures, and their comodules are entwined modules. These include  Yetter-Drinfeld and anti-Yetter-Drinfeld modules and their generalisations, hence all the modules of interest to Hopf-cyclic cohomology. In this way the interpretation of the latter as modules with flat connections [A Kaygun and M Khalkhali, Hopf modules and noncommutative differential geometry, Lett. Math. Phys. 76 (2006), 77--91] 
 is obtained as a corollary of a more general theory. We also introduce the notions of a connection in a comodule and of a bicomodule connection, and show how comodules with flat connections  can be interpreted as modules of a $C$-ring. In this way all the above mentioned Hopf modules can be interpreted as {\em comodules} with flat connections.
   \end{abstract} 
  \maketitle 

  \section{Introduction and motivation} 
The motivation for this paper comes from a recent paper \cite{KayKha:Hop} where it is proven that anti-Yetter-Drinfeld modules introduced in \cite{JarSte:coh}, \cite{HajKha:ant},  as well as $(\alpha,\beta)$-Yetter-Drinfeld modules from \cite{Pan:gen} can be understood as modules with a flat connection.  Our aim is to give an explanation of this identification in terms of corings and comodules, and to introduce a dual interpretation of anti-Yetter-Drinfeld modules as comodules with a flat connection in terms of C-rings and their modules. 

The paper is divided into two parts. The first part (sections~\ref{sec.Swe}--\ref{sec.aYD}) 
starts by describing how all the algebraic structure involved in a universal differential calculus fits in a natural way into the notion of a coring (or a coalgebra in the category of bimodules). We recall the theorem of Roiter \cite{Roi:mat} in which a bijective correspondence is established between semi-free differential graded algebras and corings with a grouplike element. A brief introduction to the theory of comodules is given and the theorem establishing a bijective correspondence between comodules of a coring with a group-like element and flat connections (with respect to the associated differential graded algebra) is given following \cite{Brz:gro}. We then specialise to corings associated to entwining structures and entwined modules.  These include all known examples of Hopf-type modules such as Hopf modules, relative Hopf modules, Long dimodules, Doi-Koppinen and alternative Doi-Koppinen modules. In particular they include Yetter-Drinfeld and anti-Yetter-Drinfeld modules and their generalisations, hence all the modules of interest to Hopf-cyclic cohomology. In this way the interpretation of the latter as modules with flat connections is obtained as a corollary of a more general theory. This part is based on the lecture delivered at the Isaac Newton Institute for Mathematical Sciences in Cambridge in August 2006, and it
does not contain any previously unknown (or unpublished) results. The aim of this part is to  put existing information together for the benefit of two mathematical communities: non-commutative geometers on one hand and algebraists working with Hopf algebras, coalgebras and corings on the other, and is intended to be easily accessible in particular to non-commutative geometers. 

The notions of an entwining structure and associated modules are formally self-dual. Rather than identifying entwined modules with comodules of a coring one can identify them with modules of an algebra in the monoidal category of bicomodules or a $C$-ring \cite[Section~6]{Brz:str}. In the second part of the paper (sections~\ref{sec.con.com}--\ref{sec.c-ring}) we first use the tools of homological coalgebra \cite{Doi:hom} to introduce the notion of a connection in a comodule and a bicomodule, and thus give a coalgebra version of topics discussed in \cite[Section~8]{CunQui:alg}. Then we interpret entwined (i.e.\ Hopf-type) modules as comodules with flat connections associated to an augmented $C$-ring. We believe that this description is new and we hope that it sheds some additional light on the structure of anti-Yetter-Drinfeld modules. We also hope that it might be useful in all situations in which $C$-rings appear as the most natural algebraic structure. One of such  areas of increasing interest is  {\em semi-infinite homological algebra} \cite{Vor:sem}, in particular semi-infinite cohomology of algebras. The basic notion of semi-regular module  \cite{Ark:sem} appearing in  this context is an example of a $C$-ring associated to an entwining structure, and, as argued in \cite{Pos:let}, $C$-rings are the natural language in which the semi-infinite cohomology of algebras should be formulated.

Throughout this paper $A$ denotes an associative unital algebra over a commutative ring $k$. Multiplication in $A$ is denoted by $\mu: A\ot A\to A$. The identity morphism for an object, say, $V$ is denoted by $V$. The unadorned tensor product is over $k$.

Recall that the {\em universal differential envelope} of $A$ is a differential graded algebra $\Omega A=\oplus_{n=0}^\infty\Omega^nA$ over $A$ (i.e.\ $A=\Omega^0A$) defined as follows. The bimodule of one-forms is 
\begin{equation}\label{one-forms}
\Omega^1 A := \ker \mu = \{\sum_i a_i\ot b_i\in A\ot A \; |\; \sum_i a_ib_i =0\}.
\end{equation}
$\Omega^1A$ has the obvious $A$-bimodule structure.
The differential  $d:A\to \Omega^1A$ is defined as
\begin{equation}\label{d}
d: a\mapsto 1\ot a- a\ot 1 = (1\ot 1)a - a(1\ot 1).
\end{equation}
One defines higher differential forms by iteration
\begin{equation}\label{higher}
\Omega^{n+1}A := \Omega^1A\otimes_A \Omega^nA,
\end{equation}
i.e.\ (or, more precisely)  $\Omega A$ is the tensor algebra of the $A$-bimodule $\Omega^1 A$, $\Omega A = T_A(\Omega^1 A)$. The differential $d$ is extended to the whole of $\Omega$ by requiring the graded Leibniz rule (and that $d\circ d =0$). This amounts to inserting the unit of the algebra $A$ in all possible places in $\Omega^nA \subset A^{\otimes n+1}$ with alternating signs.
\section{Sweedler's example and definition of corings}\label{sec.Swe}
The universal differential envelope of an algebra $A$ uses all the structure that is encoded in the notion of an algebra, i.e.\ the product (in the definition of $\Omega^1A$), the unit (in the definition of $d$) and the tensor product over $A$ (in the definition of $\Omega^nA$). In \cite{Swe:pre}, M.E.\ Sweedler proposed a different point of view on algebras. He  suggested to look at the $A$-bimodule $\cC = A\ot A$ (with the obvious $A$-actions) and consider two $A$-bilinear maps
\begin{equation}\label{swe.cop}
\DC : \cC\to \cC\ot_A\cC \simeq A\ot A\ot A, \qquad a\ot a'\mapsto a\ot 1\ot a',
\end{equation} 
\begin{equation}\label{swe.cou}
\eC :\cC\to A, \qquad \eC = \mu : a\ot a'\mapsto aa'.
\end{equation}
The algebra structure of $A$ is fully encoded in the maps \eqref{swe.cop}, \eqref{swe.cou}. It is an elementary exercise to check that the maps $\DC$ and $\eC$ make the following diagrams commute
\begin{equation}\label{coass}
\xymatrix{
 \; \cC \; \ar[rr]^{ \DC\quad} \ar[d]_{ \DC} &&\; \cC  \ot_A  
 \cC  \ar[d]^{\cC\ot_A \DC} \\  
 \cC  \ot_A\cC\ar[rr]^{\DC\ot_A\cC\quad}&&\cC\ot_A\cC\ot_A\cC\,,}
 \end{equation}
 \begin{equation}\label{counit}   
 \xymatrix{
  \; \cC \;\ar[rr]^{\DC\quad} \ar[drr]^{\simeq}  && 
  \cC \ot_A \cC  \ar[d]^{\eC\ut_A \cC}\\
    && \;A\ot_A\cC  \; , } \qquad \xymatrix{
  \; \cC  \ar[drr]^{\simeq}\ar[d]_{\DC}  & &
 \\
   \cC \ot_A \cC   \ar[rr]^{\cC\ot_A \eC} && \;\cC\ot_A A  \; . } 
   \end{equation}
Note that the diagrams \eqref{counit} simply express that $1$ is the unit in the algebra $A$. Also, note that the diagrams \eqref{coass} state that $\DC$ is a {\em coassociative map}, while \eqref{counit} state the {\em counitality axiom}. In other words, these diagrams mean that $A\ot A$ is a {\em coalgebra over a non-commutative ring $A$}. These observations lead to the following general definition (no relation of $\cC$ to $A\ot A$).
\begin{definition}\label{def.coring}
An $A$-bimodule $\cC$ is called an {\em $A$-coring} iff there are $A$-bimodule maps $\DC: \cC\to \cC\ot_A\cC$, $\eC:\cC\to A$ rendering diagrams \eqref{coass}, \eqref{counit} commutative.
\end{definition}

As for coalgebras, $\DC$ is called a {\em coproduct} and $\eC$ is called a {\em counit}. The coring $\cC=A\ot A$ is known as the {\em Sweedler} or {\em canonical} coring associated to the ring extension $k\to A$. Note in passing that $A$ itself is an $A$-coring. Thus the notion of a coring includes that of a ring.

\section{Roiter's theorem}
Going back to the universal differential envelope and realising that $\cC = A\ot A$ is a coring, we can identify $\Omega^1A$ with the kernel of the counit $\eC$. A question thus arises: are there other corings, for which the kernel of the counit gives rise to a differential graded algebra? Before this question is answered observe that  the universal differential  is defined in equation \eqref{d} as the commutator with $1\ot 1\in \cC= A\ot A$. Note that
$$
\DC(1\ot 1) = (1\ot 1)\ot_A(1\ot 1), \qquad \eC(1\ot 1) = 1.
$$
In the case of a general $A$-coring $\cC$ we can distinguish elements which have above properties and thus arrive at the following
\begin{definition}\label{def.group}
An element $g$ of an $A$-coring $\cC$ is called a {\em group-like element} provided that
$$
\DC(g) = g\ot_A g, \qquad \eC(g) =1.
$$
\end{definition}
The following remarkable result of Roiter \cite{Roi:mat} states that in fact any differential graded algebra of certain kind comes from a coring with a group-like element.
\begin{theorem}[A.V.\ Roiter]\label{thm.Roiter}
\begin{zlist}
\item Any $A$-coring $\cC$ with a group-like element $g$ gives rise to a differential graded algebra $\Omega A$ defined as follows: $\Omega^1 A = \ker\eC$,  $\Omega^{n+1}A = \Omega^1A\ot_A \Omega^nA$ and the multiplication is given by the tensor product (i.e., $\Omega A$ is the tensor algebra $\Omega A = T_A( \ker\eC)$). 
The differential
is  defined by
     $ d(a) = ga - ag$, for all $a\in A$, and, for all $c^{1}\ut_A 
\cdots \ut_A c^{n}\in (\ker\eC)^{\otimes_A n}$,
     \begin{eqnarray*}
	d(c^{1}\ut_A \cdots \ut_A c^{n}) &=&
	g\ut_A c^{1}\ut_A \cdots \ut_A c^{n}+ (-1)^{n+1}c^{1}\ut_A \cdots \ut_A c^{n}
	\ut_A g \\
	&&+ \sum_{i=1}^{n}(-1)^{i}c^{1}\ut_A  \cdots \ut_A 
	c^{i-1}\ut_A  \DC(c^{i})\ut_A  c^{i+1}
	\ut_A  \cdots\ut_A  c^{n}.
     \end{eqnarray*}
 \item A  differential graded  algebra $\Omega A$  over $A$ such that $\Omega A = T_A(\Omega^1 A)$ (that is  $\Omega^{n+1}A = \Omega^1A\ot_A \Omega^nA$; a differential graded algebra with this property is said to be {\em semi-free}), defines a coring with a grouplike element.
 \item The operations described in items (1) and (2) are mutual inverses.
 \end{zlist}
 \end{theorem}
 \begin{proof}
(1) and (3) are proven by straightforward calculations, so we only indicate how to construct a coring from a differential graded algebra (i.e.\ sketch the proof of (2)). Starting with $\Omega A$, define 
$$
\cC = Ag\oplus \Omega^1A,
$$
where $g$ is an indeterminate. In other words we define $\cC$ to be a direct sum of $A$ and $\Omega^1A$ as a left $A$-module. We now need to specify a compatible right $A$-module structure. This is defined by
$$
(ag +\omega)a' := aa'g+ada'+\omega a'.
$$
The coproduct is specified by
$$
\DC (ag) = ag\ot_A g, \qquad \DC(\omega) = g\ot_A\omega +\omega\ot_A g - d(\omega),
$$
and the counit
$$
\eC(ag+\omega) := a,
$$
for all $a\in A$ and $\omega\in \Omega^1A$. Note that this structure is chosen in such a way that $g$ becomes the required group-like element.
\end{proof}

The Roiter theorem teaches us that:\\~\\
 {\bf Semi-free differential graded algebras are in bijective correspondence with corings with a group-like element}.~\\~

The canonical coring construction can be performed for any algebra map $B\to A$ (i.e.\ it is not necessary that $B=k$) -- this is the original Sweedler's example from \cite{Swe:pre}.  In this case $\cC= A\ot_B A$ and the resulting differential graded algebra  (defined with respect to the group-like element $1\ot_B1$) corresponds to the relative universal differential forms as studied, for example, in \cite{CunQui:alg}.

\section{Comodules and flat connections}
An $A$-coring is an algebraic structure and we would like to study its (co)representations. These are given in terms of {\em comodules}. 
\begin{definition}\label{def.com}
A right $A$-module $M$ together with a right $A$-linear map $\roM: M\to M\ot_A\cC$ rendering the following diagrams 
\begin{center}
$\xymatrix{
   M\; \ar[rr]^{\varrho^M}\ar[d]_{\varrho^M}  && M\ot_A  \cC  
   \ar[d]^{M \ot_A\DC}   \\
   M\ot_A  \cC  \ar[rr]^{\varrho^M\ot_A \cC\quad} && M\ot_A\cC \ot_A\cC  
   }$\qquad 
   $\xymatrix{
   M\; \ar[rr]^{\varrho^M\quad}\ar[drr]_{\simeq} && M\ot_A\cC \ar[d]^{M\ut_A\eC} \\
    && M\ot_AA }$
\end{center}
commutative is called a {\em right $\cC$-comodule}.
\end{definition}

As for coalgebras, the map $\roM$ is called a {\em coaction}. When needed one refers to map $\roM$ which obeys the square but not the triangle  condition in Definition~\ref{def.com} as to a {\em non-counital coaction}. Comodules of the Sweedler coring $\cC = A\ot_B A$ associated to a ring extension $B\to A$ correspond bijectively to {\em descent data} for the extension $B\to A$;
see \cite[Section~25]{BrzWis:cor}. Thus corings are nowadays effectively used to describe a (generalised) noncommutative descent theory (on non-categorical level); see \cite{CaeDeG:com}.  

The existence of a group-like element in an $A$-coring $\cC$ has a very natural explanation in terms of comodules \cite{Brz:str}: $\cC$ has a group-like element if and only if $A$ is a right (or, equivalently, left) $\cC$-comodule. 

The noncommutative differential geometric  interpretation of comodules of a coring with a group-like element is provided by the following theorem taken from \cite{Brz:gro}. First recall that a {\em
    connection}  
    in a right $A$-module $M$ (with respect to a differential graded algebra $\Omega A$ over $A$) is a $k$-linear map $\nabla: 
    M\otimes_{A}\Omega^{\bullet} A\to M\otimes_{A}\Omega^{\bullet+1}A$ 
    such that, for all $\omega\in M\otimes_{A}\Omega^{k}A$ and 
    $\omega'\in\Omega A$,
    $$
    \nabla(\omega\omega') = \nabla(\omega)\omega'+(-1)^{k}\omega 
    d(\omega').
    $$
    A {\em curvature}
    of a connection $\nabla$ is a (right $A$-linear) map 
    $$
    F_{\nabla}:M\to M\otimes_{A}\Omega^{2}A,
    $$
    defined as a restriction of $\nabla\circ\nabla$ to $M$, that is,  
    $F_{\nabla}= \nabla\circ\nabla\mid_{M}$. A connection is said to 
    be {\em flat}
 if its curvature is identically equal to 0.
\begin{theorem}\label{thm.conn}
Assume that $\cC$ is an $A$-coring with a group-like element $g$, and write $\Omega A$ for the associated differential graded algebra.
\begin{zlist}
\item If $(M,\roM)$ is a right $\cC$-comodule, then the map
$$
\nabla: M\to M\ot_A\Omega^1A, \qquad m\mapsto \roM(m) - m\ot_Ag,
$$ 
is a flat connection.
\item If $M$ is a right $A$-module with a flat connection $\nabla: M\to  M\ot_A\Omega^1A$, then $M$ is a right $\cC$-comodule with the coaction
$$
\roM: M\to M\ot_A\cC, \qquad m \mapsto \nabla(m) + m\ot_A g.
$$
\item The operations described in items (1) and (2) are mutual inverses.
\end{zlist}
\end{theorem} 

This theorem is proven by a straightforward calculation and, combined with the Roiter theorem, teaches us that:
\\~\\
{\bf Flat connections with respect to a semi-free differential graded algebra are in bijective correspondence with comodules of a coring with a group-like element.}\\~

Combined with the identification of   right comodules of the Sweedler $A$-coring $A\ot_B A$ with descent data, the above observation might explain the appearance of flat connections in the descent theory cf.\ \cite{Nus:non}. In fact the correspondence between flat connections, descent theory and comodules of a Sweedler type coring goes back, at least in the commutative (algebraic geometry) case, to work of Grothendieck \cite{Gro:cry} and development of crystalline cohomology; see \cite[Chapter~2]{BerOgu:not}. Finally we would like to remark in passing that the correspondence in Theorem~\ref{thm.conn} is functorial, i.e.\ it defines an isomorphism of categories of $\cC$-comodules and $A$-modules with flat connection (with respect to $\Omega A$); see \cite[29.15--16]{BrzWis:cor} for more details.

\section{Entwined modules}\label{sec.entw}
Typically, Hopf-type modules involve data consiting of an algebra and a coalgebra, and objects which are at the same time modules and comodules with some compatibility condition. It is quite natural, therefore, to address the following problem.

Suppose that, given an algebra $A$ and a coalgebra $C$ (with coproduct $\Delta$ and counit $\eps$), we would like to construct an $A$-coring structure on $\cC = A\ot C$.  $\cC$ has an obvious left $A$-multiplication
\begin{equation}\label{psi.prod}
a(a'\ot c) := aa'\ot c, 
\end{equation}
it has also an obvious candidate for a counit, 
\begin{equation}\label{psi.cou}
\eC := A\ot \eps .
\end{equation}
In view of the identification $\cC\ot_A\cC = (A\ot C)\ot_A (A\ot C)\simeq A\ot C\ot C$, the map
\begin{equation}\label{psi.cop}
\DC := A\ot \Delta,
\end{equation}
is an obvious candidate for a coproduct for $\cC$.  To make $A\ot C$ into an $A$-coring with already specified structures \eqref{psi.prod}--\eqref{psi.cop} we need to introduce a suitable right $A$-multiplication. Obviously since $A\ot C$ must be an $A$-bimodule, in view of \eqref{psi.prod} any such a right $A$-multiplication is determined by a map $\psi: C\ot A\to A\ot C$,
\begin{equation}\label{psi}
\psi(c\ot a) := (1\ot c)a.
\end{equation}
The map $\psi$ must satisfy (four) conditions corresponding to unitality and associativity of the right $A$-multiplication and to the facts that both $\DC$ and $\eC$ are right $A$-linear maps. As observed in \cite{Brz:str} (following  a comment by M.\ Takeuchi), these four conditions are equivalent to the commutativity of the 
following  {\em bow-tie diagram}                                             
\begin{equation}\label{bow-tie}
\xymatrix{
& C\ot   A\ot   A \ar[ddl]_{\psi\otimes  A} 
\ar[dr]^{C\otimes  \mu}  &  &             
C\ot   C\ot   A \ar[ddr]^{C\ot \psi}& \\
& & C\ot   A \ar[ur]^{\Delta_C\ot  A} 
\ar[dr]^{\eps_C\ot  A}\ar[dd]^{\psi}  & &\\
A\ot  C\ot   A \ar[ddr]_{A\ot \psi} & C \ar[ur]^{C\otimes  
\iota} \ar[dr]_{\iota\otimes  C}& &A  & C\ot   A\ot   C 
\ar[ddl]^{\psi\otimes  C}\\
& &  A\ot   C \ar[ur]_{A\otimes \eps_C} \ar[dr]_{A\otimes \Delta_C} & &\\
& A\ot   A\ot   C \ar[ur]_{\mu\otimes  C} &  & A\ot  C\ot   C & , }
\end{equation}
where $\mu$ is the product in $A$ and $\iota:k\to A$ is the unit map. The map $\psi$ satisfying the conditions \eqref{bow-tie} is known as an {\em entwining map},  $C$ and $A$ are said to be {\em entwined} by $\psi$, and the triple $(A,C,\psi)$ is called an {\em entwining structure}.  These are notions introduced in this form in \cite{BrzMaj:coa} (with no reference to corings, but with an aim to recapture missing Hopf algebra symmetry needed for the construction of principal bundles over quantum homogeneous spaces). The corresponding coring $\cC = A\ot C$ is often referred to as the coring associated to an entwining structure $(A,C,\psi)$ (of course, it depends on the point of view, whether we want to see a coring as being determined by the map $\psi$ or the map $\psi$ as being determined by a coring).

One easily checks that right comodules of the $A$-coring $\cC = A\ot C$ associated to an entwining structure are simply $k$-modules $M$ which are both right $A$-modules with multiplication $\varrho_{M}: M\ot A\to M$ and right $C$-comodules with comultiplication $\roM: M\to M\ot C$ rendering commutative the following diagram 
\begin{equation}\label{ent.com}
\xymatrix{
M\otimes A \ar[rr]^{\varrho^{M}\ot A\quad} \ar[d]_{\varrho_{M}} & & 
M\otimes
C\otimes A \ar[rr]^{M\ot\psi}
& &   M\otimes A\otimes C\ar[d]^{\varrho_{M}\ot C}\\
M \ar[rrrr]^{\varrho^{M}\quad} &&& & M\otimes C \, .}
\end{equation}
Such $k$-modules are known as {\em entwined modules} (or $(A,C,\psi)$-entwined modules) and were introduced in \cite{Brz:mod}. 

Although entwining structures in this form were introduced in \cite{BrzMaj:coa}, and, at least on the first sight, the conditions expressed by the bow-tie diagram \eqref{bow-tie} might seem a bit complicated, in fact they are a special case of the structure which appeared in category theory some forty years ago and is known as  a {\em (mixed) distributive law} \cite{Bec:dis}, \cite{Van:bic}.

\section{Anti-Yetter-Drinfeld and other Hopf-type modules}\label{sec.aYD}
Since the end of the sixties, Hopf algebraists studied intensively objects with both an action and a coaction of a Hopf algebra or, more generally, with an action of an algebra and a coaction of a coalgebra which are compatible one with the other through an action/coaction of a Hopf algebra. Such objects are known as {\em Hopf-type modules}, and examples include Hopf modules of Sweedler \cite{Swe:int}, relative Hopf-modules of Doi and Takuechi \cite{Doi:str}, \cite{Tak:rel}, Doi-Koppinen Hopf modules \cite{Doi:uni},  \cite{Kop:var} or (as a special case of the latter) Yetter-Drinfeld modules \cite{RadTow:yet}, \cite{Yet:rep}.
Essentially, compatibility conditions for all known Hopf-type modules can be recast in the form of an entwining structure and are of the form of equation \eqref{ent.com}. For more information about entwining structures and their connection with Hopf-type modules we refer to \cite{CaeMil:gen} or to \cite[Section~33]{BrzWis:cor}.

The qualification {\em essentially}  appears here, since there are also variants of Hopf-type modules for {\em weak Hopf algebras} \cite{Boh:wea} (such as weak Doi-Hopf modules \cite{Boh:Doi}) and for {\em bialgebroids} \cite{Tak:gro},  \cite{Lu:alg}
(such as Doi-Koppinen modules for quantum groupoids \cite{BrzCae:Doi}). To describe the former one needs to study corings built not on $A\ot C$ but on a (left $A$-module) direct summand of $A\ot C$. Such corings are equivalently described in terms of {\em weak entwining structures} \cite{CaeDeG:mod}. To describe the latter, one works over a non-commutative ring $R$ from the onset, and studies $A$-corings on $A\ot_R C$ (to make sense of these, $C$ has to be an $R$-coring and $A$ must be an $R$-ring, i.e.\ there must be a ring map $R\to A$). These lead to {\em entwining structures over non-commutative rings} \cite{Boh:int}. In any case, to the best of author's knowledge, every known Hopf-type module (whether weak or over a non-commutative ring) is a comodule of an associated coring. This, in particular, implies to the newest additions to the family of Hopf-type modules, i.e.\ anti-Yetter-Drinfeld modules which arose naturally as coefficients in Hopf-cyclic cohomology \cite{JarSte:coh}, \cite{HajKha:ant}, and to their generalisations termed $(\alpha,\beta)$-equivariant $C$-comodules \cite{Pan:gen} \cite{KayKha:Hop}.

We illustrate the general theory of the previous sections on the example of anti-Yetter-Drinfeld modules. To this end take $A=C=H$, where $H$ is a Hopf algebra with a bijective antipode $S$. Then one can define an entwining map $\psi: H\ot H\to H\ot H$ by
\begin{equation}\label{anti.1}
\psi (c\ot a) = a\sw 2 \ot S^{-1}(a\sw 1)ca\sw 3,
\end{equation}
for all $a,c\in H$. Here $a\sw 1\ot a\sw 2 \ot a\sw 3 := (\Delta\ot H)\circ \Delta(a)$. That $\psi$ is an entwining map indeed can be easily checked by a routine calculation. While doing this exercise, the reader should notice that the only significant property (apart from multiplicativity and unitality of the coproduct) is the fact that the antipode is an anti-algebra and anti-colagebra map. Consequently, there is an $H$-coring $\cC = H\ot H$ with the right $H$-multiplication
\begin{equation}\label{anti.a}
(b\ot c) a = ba\sw 2 \ot S^{-1}(a\sw 1)ca\sw 3.
\end{equation}
The compatiblity \eqref{ent.com} for right $H$-module and $H$-comodule $M$ comes out as, for all $a\in H$,
\begin{equation}\label{anti.2}
\roM(ma) = m\sw 0a\sw 2\ot S^{-1}(a\sw 1)m\sw 1a\sw 3,
\end{equation}
where $\roM (m) = m\sw 0\ot m\sw 1$ is the $C$-coaction on $M$, i.e.\ entwined modules for \eqref{anti.1} coincide with (right-right) anti-Yetter-Drinfeld modules. Since $C=H$ is a Hopf algebra, $1_H$ is a group-like element in $H$, and hence $1_H\ot 1_H$ is a group-like element in the $H$-coring $\cC$. By the Roiter theorem there is the associated differential graded algebra and by Theorem~\ref{thm.conn} anti-Yetter-Drinfeld modules are modules with a flat connection with respect to this differential graded structure. Explicitly,
$$
\Omega^1H = \{\sum _ia_i\ot c_i\in H\ot H\; |\; \sum_i a_i\eps(c_i) =0\}.
$$
Thus, in particular $\Omega^1 H = H\ot H^+$, where $H^+ := \ker\eps$, provided $H$ is a flat $k$-module. The right $H$-action on $\Omega^1H $ is given by the formula \eqref{anti.a}. The differential comes out as
$$
d(a) = (1\ot 1) a - a(1\ot 1) = a\sw 2\ot S^{-1}(a\sw 1)a\sw 3 - a\ot 1.
$$
Note that this map is zero if $H$ is a cocommutative Hopf algebra.

Anti-Yetter-Drinfeld modules are an example of {\em ($\alpha,\beta$)-equivariant $C$-comodules} introduced in \cite{Pan:gen}. In this case $A$ is a bialgebra, $C$ is an $A$-bimodule colagebra, $\alpha: A\to A$ is a bialgebra map and $\beta: A\to A$ is an anti-bialgebra map (i.e.\ $\beta$ is both an anti-algebra and anti-coalgebra map). All these data give rise to an entwining map $\psi: C\ot A\to A\ot C$ defined by
$$
\psi (c\ot a) = a\sw 2 \ot \beta(a\sw 1)c\alpha(a\sw 3).
$$
We leave it as an exercise to work out explicitly the form of the corresponding coring $\cC = A\ot C$ and of the compatibility condition \eqref{ent.com}. If, in addition, $C$ has a group-like element $e$, then $1\ot e$ is a group-like element in $\cC$. Again, the derivation of the explicit form of the associated differential graded algebra is left as an exercise. 

\section{Connections in (bi)comodules}\label{sec.con.com}
The aim of this and the following section is to describe rudiments of the theory of connections in (bi)comodules, and to give a different interpretation of entwined modules in terms of {\em comodules} with a flat connection.  

Assume that $k$ is a field and fix a $k$-coalgebra $C$ with coproduct $\Delta: C\to C\otimes C$ and counit $\eps: C\to k$. Take a $C$-bicomodule $L$ with coactions $\Lro: L\to C\ot L$ and $\roL: L\to L\ot C$. Following \cite{Doi:hom}, a $k$-linear map $\lambda: L\to C$ is called a {\em coderivation}, provided
$$
\Delta\circ\lambda = (C\ot \lambda)\circ \Lro + (\lambda\ot C)\circ \roL.
$$
Note that $\eps\circ\lambda =0$. 

Recall that, given a right $C$-comodule $M$ with coaction $\roM$ and a left $C$-comodule $N$ with coaction $\Nro$, the cotensor product $M\square_CN$ is defined as the equaliser of $M\ot \Nro$ and $\roM\ot N$, i.e.,
$$
M\square_CN := \ker(M\ot \Nro - \roM\ot N).
$$
The assignment $(M,N)\mapsto M\square_CN$ is a functor $\M^C\times {}^C\M\to \Vect_k$. If $N$ is a $C$-bicomodule with the right coaction $\varrho^N$, then $M\square_CN$ is a right $C$-comodule with the coaction $M\square_C\varrho^N$. This is a functorial construction as well. In particular, the cotensor product makes the category of $C$-bicomodules a monoidal category. 
For any $k$-linear maps $f:M\to M'$, $g:N\to N'$ we write 
$f\square_C g$ for the restriction of  $f\ot g$ to ${M\square_C N}$. The definition of the cotensor product immediately implies that  $\roM\coten C N = M\coten C\Nro$. Note that $M\square_CC\simeq M$ and $C\square_CN\simeq N$ with isomorphisms given by the coactions and the counit. 

\begin{definition}\label{def.con}
Let $(L,\lambda)$ be a $C$-bicomodule with a coderivation. 
\begin{zlist}
\item Given a right $C$-comodule $M$, a $k$-linear map $\nabla: M\square_CL\to M$ is called a {\em connection} in $M$ with respect to $(L,\lambda)$, provided
$$
\roM\circ\nabla = (\nabla\ot C)\circ (M\square_C\roL) + M\square_C\lambda.
$$
\item Given a left $C$-comodule $N$, a $k$-linear map $\nabla: L\square_CN\to N$ is called a {\em connection} in $N$ with respect to $(L,\lambda)$, provided
$$
\Nro\circ\nabla = (C\ot \nabla)\circ (\Lro\square_CN) + \lambda\square_CN.
$$
\end{zlist}
\end{definition}
With any coalgebra $C$ one can associate the {\em universal coderivation}. Write $C^+ := \ker \eps$, and set $L(C):= C\ot C^+$. View $L(C)$ as a left $C$-comodule via $\Delta\ot C^+$ and as a right $C$-comodule with the coaction $C\ot\Delta|_{C^+} - \Delta\ot C^+$. Then the map
$$
\lambda_C : L(C)\to C, \qquad c\ot d\mapsto \eps(c)d,
$$
is a coderivation, called a {\em universal coderivation}. Note that equivalently $L(C)$ can be defined as $\coker\, \Delta$, with $\lambda : \overline{c\ot d}\mapsto \eps(c)d - c\eps(d)$, where  $\bar{x}$ denotes the element of $\coker\, \Delta$ corresponding to $x\in C\ot C$.  The isomorphism $\coker\,\Delta\to C\ot C^+$ is given by $\overline{c\ot d}\mapsto c\ot d - \Delta(c)\eps(d)$ (with the inverse $c\ot d\mapsto \overline{c\ot d}$). This identification of $L(C)$ with $\coker\, \Delta$ makes the duality between $L(C)$ and $\Omega^1A$ in \eqref{one-forms} more transparent. 
\begin{lemma}\label{lem.uni}
A right (resp.\ left) $C$-comodule  admits a connection with respect to $(L(C),\lambda_C)$ if and only if it is an injective $C$-comodule.
\end{lemma}
\begin{proof}
This is a dual version of the characterisation of modules with a connection with respect to universal differential structure; see \cite[Corollary~8.2]{CunQui:alg}. Write $\pi: C\ot C\to \coker\, \Delta$ for the canonical epimorphism. Given a connection $\nabla: M\square_CL(C)\simeq M\square_C \coker\,\Delta \to M$ in a right $C$-comodule $M$, the map
$$
\sigma_r: M\ot C\to M, \qquad \sigma_r:= M\ot\eps + \nabla\circ(M\ot\pi)\circ(\roM\ot C),
$$
is a right $C$-colinear retraction of the coaction $\roM$.
Hence $M$ is an injective comodule. Conversely, if $M$ is an injective $C$-comodule and $\sigma_r$ is a right $C$-colinear retraction of $\roM$, then $\nabla:= \sigma_r\circ (M\coten C(\eps\ot C^+))$  is a connection in $M$. 

If $\nabla : L(C)\square_CN\simeq  (\coker\,\Delta)\square_C N \to N$ is a connection in a left $C$-comodule $N$, the map
$$
\sigma_l: C\ot N\to N, \qquad \sigma_l:= \eps\ot N - \nabla\circ(\pi\ot N)\circ(C\ot\Nro),
$$
is a left $C$-colinear retraction of the coaction $\Nro$.
Hence $N$ is an injective comodule. Conversely, if $N$ is an injective $C$-comodule and $\sigma_l$ is a left $C$-colinear retraction of $\Nro$, then $\nabla:= \sigma_l\circ((\eps\ot C^+)\coten C N)$  is a connection in $N$.  
\end{proof}

If $M$ is a $C$-bicomodule, one can consider connections in $M$ as a right and left $C$-comodule, and demand compatibility with other comodule structures. Dualising definitions in \cite[Section~8]{CunQui:alg}, we can thus propose
\begin{definition}\label{def.con.bic}
Let $(L,\lambda)$ be a $C$-bicomodule with a coderivation, and let $M$ be a $C$-bicomodule.
\begin{zlist}
\item  A left $C$-colinear connection in a right $C$-comodule $M$ (with respect to $(L,\lambda)$) is called a {\em right connection} in $M$ with respect to $(L,\lambda)$.
\item A right $C$-colinear connection in a left $C$-comodule $M$ (with respect to $(L,\lambda)$) is called a {\em left connection} in $M$ with respect to $(L,\lambda)$.
\item A {\em bicomodule connection} in $M$ is a pair $(\nabla_l,\nabla_r)$ such that $\nabla_r$ is a right connection in $M$ and $\nabla_l$ is a left connection in $M$.
\end{zlist}
\end{definition} 

Similarly to Lemma~\ref{lem.uni}, the existence of bicomodule connections with respect to the universal coderivation is closely related to injectivity. A $C$-bicomodule $M$ is {injective} if and only if there exists a $C$-bicolinear retraction of
$$
(\Mro\ot C)\circ\roM = (C\ot \roM)\circ \Mro,
$$
where $\roM: M\to M\ot C$ and $\Mro: M\to C\ot M$ are coactions. Thus if $M$ is an injective $C$-bicomodule, it is also injective as a left and right $C$-comodule. Furthermore, there exist $C$-bicolinear retractions of coactions. 

\begin{proposition}\label{prop.inj.bic}
A $C$-bicomodule $M$ admits a bicomodule connection with respect to the universal coderivation $(L(C),\lambda_C)$ if and only if $M$ is an injective bicomodule.
\end{proposition}
\begin{proof} 
This is dual to \cite[Proposition~8.3]{CunQui:alg}. Let $(\nabla_l, \nabla_r)$ be a bicomodule connection. Since a right connection  $\nabla_r$ is left colinear, the corresponding retraction $\sigma_r$ as constructed in the proof of Lemma~\ref{lem.uni} is also left $C$-colinear. Consequently, $\sigma_r\circ(\sigma_l\ot C)$, where $\sigma_l$ is a retraction of $\Mro$ corresponding to $\nabla_l$, is a $C$-bicolinear retraction of $(\Mro\ot C)\circ\roM$. The converse follows immediately by Lemma~\ref{lem.uni} and the discussion after Definition~\ref{def.con.bic}. 
\end{proof}

In particular, $L(C)$ admits a bicomodule connection with respect to coderivation $(L(C),\lambda_C)$ if and only if $C$ is a formally smooth coalgebra; see \cite[Theorem~1.4]{JarLle:her}.

To define a torsion and curvature of a connection, we need to consider extended coderivations.

\begin{definition}\label{def.ext}
Let $(L,\lambda)$ be a $C$-bicomodule with a coderivation.
By an {\em extended coderivation} we mean a triple $(L,\lambda,\lambda')$, where   $\lambda': L\square_C L\to L$ is  a $k$-linear map such that $\lambda'\circ\lambda =0$, $\lambda'$ is a connection in the left $C$-comodule $L$ with respect to $(L,\lambda)$, and a connection in the right  $C$-comodule $L$ with respect to $(L,-\lambda)$. Explicitly, we require 
\begin{equation}\label{exten}
\Lro \circ \lambda' = (C\otimes \lambda')\circ (\Lro\square_C L)+ \lambda\square_C L.
\end{equation}
\begin{equation}\label{exten.r}
\roL \circ \lambda' = (\lambda'\otimes C)\circ (L\square_C\roL )- L\square_C \lambda.
\end{equation}
\end{definition}
An extended coderivation gives rise to a chain complex
$$
\xymatrix{
L\square_C L \ar[r]^{\lambda'} & L\ar[r]^\lambda & C\ar[r]^\eps & k.}
$$
 In order not to clatter the notation $\lambda'$ is simply denoted by $\lambda$, and we write $(L,\lambda)$ for $(L,\lambda,\lambda')$.  The universal coderivation can be extended, for $L(C)\square_C L(C)\simeq C\ot C^+\ot C^+$ and the extension of $\lambda_C$ can be defined by  $\lambda_{C} : c\ot c'\ot c''\mapsto \eps(c)c'\ot c''$. (While checking \eqref{exten} and \eqref{exten.r} the reader should note that one needs to view $c\ot c'\ot c''$ in $L(C)\coten C L(C)$ using the right $C$-comodule structure of $C\ot C^+$.)
 
 Since $L$ is itself a $C$-bicomodule, one can study bicomodule connections in $L$. In case $(L,\lambda)$ is an extended coderivation, this becomes very simple.
 
 \begin{proposition} \label{prop.l-r}
 Let $(L,\lambda)$ be an extended coderivation. Then the formula
 $$
 \nabla_l = \lambda +\nabla_r.
 $$
 gives a bijective correspondence between left and right connections in the $C$-bicomodule $L$ with respect to $(L,\lambda)$.
 \end{proposition}
 \begin{proof} 
 This is a bicomodule version of \cite[Proposition~8.5]{CunQui:alg}. Take a right connection $\nabla_r$ in $L$ and set $\nabla_l = \lambda +\nabla_r$. Then
\begin{eqnarray*}
 \Lro\circ\nabla_l 
 &=& (C\otimes \lambda)\circ (\Lro\square_C L)+ \lambda\square_C L +(C\otimes \nabla_r)\circ (\Lro\square_C L)\\
 &=& (C\otimes \nabla_l)\circ (\Lro\square_C L)+ \lambda\square_C L,
\end{eqnarray*}
 where the first equality follows by \eqref{exten} and by the left $C$-colinearity of $\nabla_r$. Hence $\nabla_l$ is a connection in the left $C$-comodule $L$. Furthermore,
 \begin{eqnarray*}
 \roL\circ\nabla_l 
 &=& (\lambda\otimes C)\circ (L\square_C\roL )- L\square_C \lambda + 
 (\nabla_r\ot C)\circ (L\square_C\roL) + L\square_C\lambda\\
 &=& (\nabla_l\ot C)\circ (L\coten C \roL),
 \end{eqnarray*}
 where the first equality follows by \eqref{exten.r} and by the definition of a connection in a right $C$-comodule. Thus $\nabla_l$ is a left connection in $L$. The fact that a left connection induces the right connection is proven in a similar way. The bijectivity is obvious.
\end{proof}

Using the same arguments as in the proof of Proposition~\ref{prop.l-r}, one easily checks that if $(L,\lambda)$ is an extended derivative and $(\nabla_l,\nabla_r)$ is a bicomodule connection in $L$ (with respect to $(L,\lambda)$), then the map
$$
T_{(\nabla_l,\nabla_r)} : L\coten CL \to L, \qquad T_{(\nabla_l,\nabla_r)} := \nabla_l - \lambda -\nabla_r,
$$
is a $C$-bicomodule map. $T_{(\nabla_l,\nabla_r)}$ is called a {\em torsion} of $(\nabla_l,\nabla_r)$. 
Proposition~\ref{prop.l-r} implies that, given an extended coderivation $(L,\lambda)$ any left (or right) connection in $L$ gives rise to a {\em torsion-free} bicomodule connection. In particular this is true for the universal (extended) coderivation $(L(C),\lambda_C)$. Thus, in view of Proposition~\ref{prop.inj.bic},  to prove that $L(C) = C\ot C^+$ is an injective bicomodule suffice it to find a bicolinear retraction of one of the $C$-coactions in $L(C)$. 

Recall that a coalgebra $C$ is said to be {\em coseparable} if there exists a $k$-linear map $\delta :C\ot C\to k$, such that
$$
(\delta \ot C)\circ (C\ot \Delta) = (C\ot \delta)\circ (\Delta \ot C), \qquad \delta\circ\Delta = \eps.
$$
Such a map $\delta$ is called a {\em cointegral}. Every bicomodule of a coseparable coalgebra is injective, 
so any $C$ bicomodule has a bicomodule connection with respect to $(L(C),\lambda_C)$. In particular, if $C$ is a coseparable coalgebra, then a right connection in $L(C)$ with respect to $(L(C),\lambda_C)$ can be defined as
$$
\nabla_r = (C\ot C\ot \delta)\circ (C\ot \Delta\ot C^+) - \Delta\ot\delta \, .
$$
The corresponding left connection computed from Proposition~\ref{prop.l-r}  is
$$
\nabla_l = \eps\ot C^+\ot C^+ + (C\ot C\ot \delta)\circ (C\ot \Delta\ot C^+) - \Delta\ot\delta.
$$
The resulting bicomodule connection is torsion-free (i.e.\ $T_{(\nabla_l,\nabla_r)}=0$).

\begin{lemma}\label{lem.exten}
\begin{zlist}
\item Let $\nabla$ be a connection in a right $C$-comodule $M$ with respect to an extended coderivation $(L,\lambda)$. Define
$$
\nabla_\lambda: M\square_C L\square_C L\to M\ot L, \qquad \nabla_\lambda:= \nabla\square_C L + M\square_C \lambda.
$$
Then $\im \nabla_\lambda \subseteq M\square_C L$.
\item Let $\nabla$ be a connection in a left $C$-comodule $N$ with respect to an extended coderivation $(L,\lambda)$. Define
$$
\nabla_\lambda: L\square_C L\square_C N\to L\ot N, \qquad \nabla_\lambda:=  \lambda\square_C N -L\square_C \nabla.
$$
Then $\im \nabla_\lambda \subseteq L\square_C N$.
\end{zlist}
\end{lemma}
\begin{proof}
This is proven by a straightforward calculation which uses the definitions of the cotensor product and connection, and equations \eqref{exten}, \eqref{exten.r}.
\end{proof} 

In view of Lemma~\ref{lem.exten} it is possible to make the following
\begin{definition}\label{def.flat}
Let $(L,\lambda)$ be an extended coderivation. The {\em curvature} of a connection $\nabla$ in a right (resp.\ left) $C$-comodule is defined as 
$$
 F_\nabla:= \nabla\circ\nabla_\lambda.
$$
The connection $\nabla$ is said to be {\em flat} if its curvature vanishes, $F_\nabla=0$.
\end{definition}
Any flat connection in a right $C$-comodule $M$ (resp.\ left $C$-comodule $N$) gives rise to  a chain complex
$$
\xymatrix{
M\square_CL\square_C L \ar[r]^{\nabla_\lambda} & M\square_CL\ar[r]^\nabla & M,} \qquad \mbox{(resp.}\ \xymatrix{L\square_CL\square_C N \ar[r]^{\nabla_\lambda} & L\square_CN\ar[r]^\nabla & N}).
$$

\section{Modules of $C$-rings and flat connections} \label{sec.c-ring}
In Section~\ref{sec.entw}, $(A,C,\psi)$-entwined modules were identified with comodules of an associated coring. As observed in \cite[Proposition~6.2]{Brz:str}, equivalently, one can describe $(A,C,\psi)$-entwined modules as modules of the $C$-ring or the monoid in the category of $C$-bicomodules associated to $(A,C,\psi)$. The aim of this section is to show that any $C$-ring with a character gives rise to an extended coderivation (the dual Roiter theorem), and that the modules of this $C$-ring can be identified with comodules with a flat connection. This, in particular, gives an interpretation of entwined modules (in case, when $A$ has a character), hence anti-Yetter-Drinfeld and Yetter-Drinfeld modules, as comodules with flat connections.  

We assume that $k$ is a field and $C$ is a $k$-coalgebra. Let $\cA$ be a $C$-bicomodule with coactions $\Aro: \cA\to C\ot \cA$ and $\roA:\cA\to \cA\ot C$. $\cA$ is called a {\em $C$-ring} if there are  two bicomodule maps $\mA:
\cA\coten C\cA\to\cA$ and $\eA:C\to \cA$ such that
$$
\xymatrix{\cA\coten C\cA\coten C \cA \ar[rr]^{\mA\coten C\cA}\ar[d]_{\cA\coten C\mA} && \cA\coten C \cA\ar[d]^\mA\\
\cA\coten C \cA\ar[rr]^\mA && \cA \, ,}
$$
$$
\xymatrix{\cA \ar[rr]^{\roA}\ar[d]_\cA && \cA\coten C C \ar[d]^{\cA\coten C \eA} & & \cA \ar[rr]^{\Aro}\ar[d]_\cA && C\coten C \cA \ar[d]^{\eA\coten C\cA} \\
\cA && \cA\coten C \cA \ar[ll]_\mA \, , && \cA && \cA\coten C \cA \ar[ll]_\mA \, .}
$$
A {\em character} in $\cA$ is a $k$-linear map $\kappa: \cA\to k$ such that
\begin{equation}\label{char}
\kappa\circ\mA = \kappa\coten C\kappa, \qquad \kappa\circ\eA = \eps.
\end{equation}
 A right module of a $C$-ring $\cA$
is a right $C$-comodule $M$ with coaction $\roM: M\to M\ot C$ together with a right $C$-comodule map
$\wra: M\coten C \cA\to M$ such that
\begin{equation}\label{action}
\xymatrix{M\coten C\cA\coten C \cA \ar[rr]^{\wra\coten C\cA}\ar[d]_{M\coten C\mA} && M\coten C \cA\ar[d]^\wra && M \ar[rr]^{\roM}\ar[d]_M && M\coten C C \ar[d]^{M\coten C \eA} \\
M\coten C \cA\ar[rr]^\wra && M \, , && 
M && M\coten C \cA \ar[ll]_\wra \, .}
\end{equation}
The map $\wra$ is called a {\em right $\cA$-action}. $C$ is a right $\cA$-module if and only if there is a character in $\cA$. 

A $C$-ring can be understood as an associative unital algebra (a monoid) in a monoidal category of $C$-bicomodules. In view of this, the first assertion of the following proposition can be viewed as a standard result in the theory of algebras.
\begin{proposition}\label{prop.complex}
Let $\cA$ be a $C$-ring with a character $\kappa$. Set $\bA := \coker\, \eA$ and consider the following (left infinite) sequence
\begin{equation}\label{seq}
\xymatrix{
...\ar[r]^{\lambda\hspace{.5cm}} & \bA \coten C \bA\coten C\bA \ar[rr]^\lambda && \bA \coten C \bA\ar[r]^\lambda & \bA\ar[r]^\lambda & C\ar[r]^\eps &k,
}
\end{equation}
with the maps $\lambda$ defined as follows. Write $\pi: \cA\to \bA$ for the canonical surjection. Then $\lambda: \bA\to C$ is defined by 
$$
\lambda\circ\pi = (\kappa\ot C)\circ\roA - (C\ot\kappa)\circ\Aro,
$$
and $\lambda: \bA^{\square_C n}\to \bA^{\square_C n-1}$,
$$
\lambda\circ\pi^{\square_C n} = \kappa\square_C\pi^{\square_C n-1} + \sum_{l=1}^{n-1} (-1)^l\pi^{\square_C l-1}\square_C(\pi\circ\mA)\square_C \pi^{\square_C n-l-1}+ (-1)^n \pi^{\square_C n-1}\square_C \kappa.
$$
Then \eqref{seq} is a chain complex and $(\bA,\lambda)$ is an extended coderivation. 
\end{proposition}
(Compare the definition of $\lambda$ with that of $d$ in Theorem~\ref{thm.Roiter}.)

\begin{proof}
Note that the map $\pi$ has a left $C$-comodule section $\sigma: \bA\to \cA$ given by
$$
\sigma\circ\pi = \cA - (\eA\ot \kappa)\circ\Aro.
$$
This implies that, for all right $C$-comodules $M$, $M\coten C \bA = \im\, (M\coten C\pi)$. In particular $\bA^{\square_C n} = \im\,\left( \pi^{\square_C n}\right)$, and the definitions of the $\lambda$ are justified.  We only check that $(\bA,\lambda)$ is an extended coderivation. The $C$-coactions $\robA$ and $\bAro$ on $\bA$ are induced from $\roA$ and $\Aro$ by $\pi$, hence the map $\pi$ is $C$-bicolinear, i.e., $\robA\circ \pi = (\pi\ot C)\circ\roA$ and $\bAro\circ\pi = (C\ot \pi)\circ\Aro$. In view of this one can compute
\begin{eqnarray*}
[(C\ot\lambda)\circ (\bAro\coten C \bA)+ \lambda\coten C \bA]\circ (\pi\coten C\pi) \!\!\!&=& \!\!\! (C\ot \kappa \coten C \pi - C\ot \pi\circ\mA + C\ot \pi\coten C\kappa ) \circ (\Aro\coten C \cA)\\
&&+ (\kappa\ot C\ot \cA)\circ (\roA\coten C \pi) - (C\ot\kappa\ot \cA)\circ (\Aro\coten C \pi)\\
&=& \!\!\! \bAro\circ (\pi\coten C\kappa - \pi\circ\mA + \kappa\coten C\pi)
= \bAro\circ\lambda
\circ (\pi\coten C\pi),
\end{eqnarray*}
where the first equality follows by the $C$-colinearity of $\pi$ and definitions of the $\lambda$, while the second follows by the $C$-colinearity of the multiplication $\mA$ and $\pi$, and by the definition of the cotensor product. This proves the condition \eqref{exten}. Equation \eqref{exten.r} is proven by similar arguments.
\end{proof}

As a particular example of construction in Proposition~\ref{prop.complex} one can derive the complex associated to the universal coderivation. Simply view $C\otimes C$ as a $C$-ring with multiplication 
$$C\otimes C\square_C C\otimes C\simeq C\otimes C\otimes C\to C\otimes C, \qquad c\otimes c'\otimes c''\mapsto \eps(c') c\otimes c'',
$$
unit $\Delta$ and character $\eps\otimes\eps$. Identifying $\overline{C\ot C} = \coker\,\Delta$ with $L(C) = C\ot C^+$ we obtain
$$
\lambda: L(C)^{\square_Cn} \simeq C\ot (C^+)^{\ot n}\to C\ot (C^+)^{\ot n-1}\simeq L(C)^{\square_Cn-1}, \;\;\; c_0\ot\cdots \ot c_n\mapsto \eps(c_0)c_1\ot\cdots \ot c_n.
$$

\begin{theorem}\label{thm.flat}
Let $\cA$ be a $C$-ring with a character $\kappa$ and let $(\bA,\lambda)$, $\pi$  be as in Proposition~\ref{prop.complex}. For any right $C$-comodule $M$, the formula
$$
\wra = \nabla\circ (M\square_C\pi) + M\coten C\kappa,
$$
gives a bijective correspondence between right $\cA$-actions on $M$ and flat connections in $M$ with respect to $(\bA,\lambda)$.
\end{theorem}
\begin{proof} Similarly to Theorem~\ref{thm.conn} this is proven by direct calculations. For example, suppose $M$ is a right $\cA$-module with $C$-coaction $\roM$ and $\cA$-action $\wra$, then 
$$
(\wra-M\coten C \kappa)\circ(M\coten C\eA)\circ \roM = M - (M\ot\eps)\circ\roM =0,
$$
by the second of diagrams \eqref{action}, the second of equations \eqref{char} and the counitality of the coaction. This means that the there is a unique map $\nabla$ such that $\nabla\circ (M\square_C\pi) =\wra -M\coten C \kappa$. Next, 
\begin{eqnarray*}
\roM\circ \nabla \circ(M\coten C\pi) &=& \roM\circ\wra - \roM\circ (M\coten C \kappa)\\
&=& (\wra\ot C)\circ (M\coten C\roA) - (M\ot C\ot\kappa )\circ (M\coten C \Aro)\\
&=& [(\wra- M\coten C \kappa)\ot C]\circ (M\coten C\roA)
  + M\coten C [ (\kappa\ot C)\circ \roA- (C\ot\kappa )\circ \Aro]\\
&=& [(\nabla\ot C)\circ (M\coten C\robA) + M\coten C \lambda] \circ (M\coten C\pi) .
\end{eqnarray*}
where the second equality follows by the right $C$-colinearity of the action $\wra$, and by the definition of the cotensor product. The final equality is a consequence of the fact that, by construction, $\pi$ is a right $C$-comodule map. The flatness of $\nabla$ is a straightforward consequence of the associativity of $\wra$, the definition of $\lambda$ and equations \eqref{char}. The verification that given a flat connection the formula in the theorem gives an $\cA$-action is left to the reader.
\end{proof}

Starting with an entwining structure $(A,C,\psi)$, one constructs a $C$-ring $\cA = C\ot A$ with the $C$-coactions 
$$
\Aro = \Delta\ot A, \qquad \roA = (C\otimes\psi)\circ (\Delta\ot A),
$$
multiplication 
$$
\mA: (C\otimes A)\square_C (C\otimes A) \simeq C\otimes A\otimes A\to C\otimes A, \qquad c\otimes a\otimes a'\mapsto c\otimes aa',
$$
and unit $\eA: C\to C\otimes A$, $c\mapsto c\otimes 1$. The category of entwined modules is then isomorphic to the category of right $\cA$-modules; see \cite[Proposition~6.2]{Brz:str}. If $\chi: A\to k$ is a character, then $\kappa:= \eps\otimes \chi$ is a character in the $C$-ring $\cA$. The cokernel $\bA$ of $\eA$ can be identified with $C\otimes \bar{A}$, where $\bar{A}= A/k1$. Thus $\bA\coten C\bA \simeq C\ot\bar{A}\ot\bar{A}$. With this identification and writing $\ba$ for the image of $a\in A$ under the canonical surjection $A\to \bar{A}$, the coderivative and its extension come out as, for all $a,b\in A$, $c\in C$,
$$
\lambda(c\ot \ba) = \sum_\psi \chi(a_\psi)c^\psi -\chi(a)c, \quad
\lambda(c\ot \ba \ot \bar{b}) = \sum_\psi \chi(a_\psi)c^\psi\ot \bar{b} - c\ot \overline{ab} + c\ot\ba\chi(b),
$$
where $\psi(c\ot a) = \sum_\psi a_\psi\ot c^\psi$.

In the case of anti-Yetter-Drinfeld modules of a Hopf algebra $H$, the corresponding $H$-ring is $\cA_{aYD}=H\otimes H$, hence there is a character $\kappa: =\eps\ot\eps$.  $\bA_{aYD}=H\otimes \bar{H}$, and, in view of  \eqref{anti.1}, the coderivative and its extension come out as, for all $a,b,c\in H$, 
$$
\lambda_{aYD}(c\ot \ba) = S^{-1}(a\sw 1)ca\sw 2 -\eps(a)c, \;\;\;
\lambda_{aYD}(c\ot \ba \ot \bar{b}) = S^{-1}(a\sw 1)ca\sw 2\ot \bar{b} - c\ot \overline{ab} + c\ot\ba\eps(b).
$$
In particular $\lambda_{aYD}: H\otimes\bar{H}\to H$ vanishes if $H$ is a commutative Hopf algebra. For any anti-Yetter-Drinfeld module $M$, the corresponding flat connection in $H$-comodule $M$ is, for all $a\in H$ and $m\in M$,
$$
\nabla(m\ot\ba) = ma - m\eps(a).
$$

\section{Comments on semi-group-like elements and conventions.}

An element $g$ of an $A$-coring $\cC$ is called a {\em semi-group-like element} provided $\DC(g) = g\ot_A g$.  For any $A$-coring $\cC$ with a semi-grouplike-element $g$ one can construct a (semi-free) differential graded algebra over $A$, by setting $\Omega^1 A =\cC$, $\Omega A = T_A(\cC)$ with the same formulae for $d$ as in  Theorem~\ref{thm.Roiter}(1) (cf.\ \cite[29.2]{BrzWis:cor}).\footnote{In view of Theorem~\ref{thm.Roiter}(2), the differential graded algebra $\Omega A$ can be understood as an algebra coming from a coring with a group-like element (this coring is then obtained as a direct sum of $\cC$ with a left $A$-module $A$).} 
The same formula as in Theorem~\ref{thm.conn}(1) assigns a flat connection $\nabla: M\to M\ot_A\cC$ to a right coaction $\roM$. Both formulae (1) and (2) in Theorem~\ref{thm.conn} establish a bijective correspondence between {\em non-counital coactions} and flat connections with respect to the differential graded algebra $\Omega A = T_A(\cC)$. Specialising to (anti-)Yetter-Drinfeld modules or $(\alpha, \beta)$-equivariant $C$-comodules one then obtains the results of \cite{KayKha:Hop}.

Throughout this paper we prevalently used the {\em right-right conventions}, i.e.\ we studied right actions and right coactions. Obviously, one can study left comodules over an $A$-coring (these will correspond to left $A$-modules with a flat connection) or left modules over a $C$-ring (these will correspond to left $C$-comodules with a flat connection). In the case of entwining structures, there are four possible conventions (right-right, right-left, left-right, left-left); thus, for example, there are four types of entwining structures corresponding to four types of anti-Yetter-Drinfeld modules. One can move freely between these conventions by using opposite/co-opposite algebras and/or coalgebras; see \cite{CaeMil:gen}. Obviously, although this requires some care, it does not introduce any new non-trivial features. 

\section*{Acknowledgements}
I am grateful to the organisers of the Noncommutative Geometry Programme at the Isaac Newton Institute for the Mathematical Sciences, Cambridge, and 
the organisers of the conference on  ``New techniques in Hopf algebras and graded ring theory", Brussels, for the invitations and financial support. I would like to thank Zoran \v Skoda for a discussion on crystalline cohomology and Dmitriy Rumynin for bringing  semi-infinite cohomology to my attention.

  \end{document}